\documentstyle{amsppt}
\input epsf
\magnification=\magstep1
\vsize=9.0 true in
\hsize=6.5 true in
\topmatter
\title
Mean values of $L$--functions and symmetry
\endtitle
\rightheadtext{Mean values and symmetry}
\author
J.B.~Conrey \\
D.W.~Farmer 
\endauthor
\thanks Research supported 
by the American Institute of Mathematics.
Research of the first author supported in
part by a grant from the NSF.
\endthanks
\address
American Institute of Mathematics,
Palo Alto, CA 94306
\endaddress
\abstract
Recently Katz and Sarnak introduced the idea of a symmetry
group attached to a family of $L$--functions, and they 
gave strong  evidence that the symmetry group governs many
properties of the distribution of zeros of the $L$--functions.
We consider the mean--values of the $L$--functions and the 
mollified mean--square of the $L$--functions and
find evidence that these are also governed by the symmetry group.
We use recent work of Keating and Snaith to 
give a complete description of these mean values.
We find a connection to the
Barnes--Vign\'eras $\Gamma_2$--function and to 
a family of self--similar functions.
\endabstract
\endtopmatter
\document

\NoBlackBoxes

\def\mod{\bmod}

\def\R {{\Bbb R}}

\def\intl{\int\limits}
\def\thalf{{\textstyle{1\over 2}}}

\def\({\left(}
\def\){\right)}
\def\[{\left[}
\def\]{\right]}

\head 
1. Introduction
\endhead

Katz and Sarnak \cite{KS} have introduced 
the idea of a family of $L$--functions with an associated
symmetry type.  The symmetry type has been shown to
govern the distribution and spacing of zeros in
the function field case \cite{KS2}, and there
is strong numerical evidence \cite{Ru} that it
governs the behavior of zeros in more general cases.
The paper \cite{ILS} also shows clear evidence that the symmetry type
governs the distribution of low--lying zeros
of a wide variety of $L$--functions.
In this paper we give evidence that the symmetry type of a 
family of $L$--functions governs the behavior of
mean values of the $L$--functions.  

The most well--understood mean values are the $2k$th moments of the 
Riemann $\zeta$--function.  We will describe that case in detail,
and then compare this to the examples provided by the
families of Katz and Sarnak.

It is a folklore conjecture
that for every $k\ge 0$ there is a constant $c_k$
such that
$$
I_k(T):=
\frac 1T \intl_0^T\left|\zeta(\textstyle{1\over2}+it)\right|^{2k}~dt
\sim c_k \log^{k^2}T
$$
as $T\to \infty$. 
Classically, Hardy and Littlewood \cite{HL} proved that
$$
I_1(T)\sim \log T,
$$
and Ingham \cite{I} showed
$$
I_2(T)\sim \frac{1}{2\pi^2}\log^4 T,
$$
so $c_0=c_1=1$ and $c_2=1/2\pi^2.$
No other mean values of the $\zeta$-function have been established.

More information can be given about $c_k$, for Conrey and Ghosh \cite{CG1} 
showed (assuming the Riemann Hypothesis) that for all $k\ge 0$,
$$
I_k \ge (1+o(1))\frac{a_k}{\Gamma(1+k^2)} \log^{k^2}T,
$$
where $a_k$ is an arithmetic constant given by
$$
a_k=\prod_p(1-1/p)^{k^2}\sum_{j=0}^\infty \frac{ d_k(p^j)^2}{p^j}.
\eqno (1.1)
$$
Here $d_k(n)$ is the $n$th coefficient in the Dirichlet series for $\zeta(s)^k$; it is multiplicative and is given by
$$
d_k(p^j)=\frac{\Gamma(k+j)}{\Gamma(k)j!}.
$$
Gonek \cite{G} extended the above result to include 
all $k>-{1\over 2}$.

Conrey and Ghosh \cite{CG3} established properties of~$a_k$. 
These included that $a_k$, regarded as a function of~$k$, 
is entire of order 2, satisfies a symmetry $a_k=a_{1-k}$, 
and has all of its zeros on the line $\Re k =\frac12$. 
In addition, all of the zeros of
the derivative $\frac{d}{dk} a_k$  have real parts 
equal to $\frac12$ with two exceptions $a_0'=a_1'=0$.
Conrey and Gonek \cite{CGo} give an asymptotic formula for   
$\log a_k$ as $k \to \infty$. 
This formula is relevant to understanding the 
extremely large values of $\zeta(s)$.

Conrey and Ghosh \cite{CG3} defined the number $g_k$  
implicitly by the (conjectural) formula
$$
I_k(T)\sim g_k \frac{a_k}{\Gamma(1+k^2)} \log^{k^2}T.
\eqno (1.2)
$$ 
The quantity $g_k$ is natural in the following sense.
All approaches to mean value theorems for Dirichlet
series have relied on techniques from the theory of 
Dirichlet polynomials. The main tool there is the 
mean value theorem of Montgomery and Vaughan \cite{MV}:
$$
\intl_0^T\left|\sum_{n=1}^N a_n n^{it}\right|^2~dt=
\sum_{n=1}^N (T+O(n))|a_n|^2.
$$
To use this  to obtain asymptotic formulae it is
usually necessary to have $N\ll T$. 
Thus, it is natural in some sense to measure the mean square of
$\zeta(s)^k$ against the mean value of
a Dirichlet polynomial approximation to $\zeta(s)^k$
using a polynomial of length $T$.
It is easy to show that
$$
\frac 1T \intl_0^T \left|\sum_{n \le T}\frac{d_k(n)}{n^{1/2+it}}\right|^2~dt
\sim \frac{a_k}{\Gamma(1+k^2)} \log^{k^2}T.
$$
Thus, $g_k$ is a measure of ``how many polynomials 
of length $T$ are needed to capture the mean square
of $\zeta(s)^k$.''

The classical results of Hardy and Ingham can be phrased as $g_1=1$ and $g_2=2$. Recently, Conrey and 
Ghosh \cite{CG2} made the conjecture that $g_3=42$. Still more recently,
Conrey and Gonek \cite{CGo} conjectured that $g_4=24024$.
The conjectures for $g_3$ and $g_4$ are based on
Dirichlet polynomial techniques.
With current methods
it is probably not possible to use those techniques to
conjecture $g_k$ for larger values of $k$.
We will return to the function $g_k$ after 
giving a general discussion of mean
values and a description of the 
situation for the families of $L$--functions of Katz and Sarnak.

We will now describe the situation
for some of the families of Katz and Sarnak. 
The mean values will be of the shape
$$
\frac 1{{\Cal Q}^*} \sum_{f\in \Cal F \atop c(f)\le\Cal Q } 
V(L_f({\textstyle \frac 12}))^k\sim
 g_k \frac{a_k}{\Gamma(1+B(k))} (\log {\Cal Q}^A)^{B(k)}.
\eqno (1.3)
$$
The $L$--functions are normalized to have a functional
equation $s\leftrightarrow 1-s$, so $L_f(\frac 12)$ is the ``critical value.''
Here we think of the family $\Cal F$ as being partially ordered
by ``conductor'' $c(f)$,  
with ${\Cal Q}^*$ the number of elements with 
$c(f)\le\Cal Q$.  
We set $V(z)=z$ or~$|z|^2$ depending on
the symmetry type of the family.
The parameters $g_k$ and $B(k)$ depend only on  the symmetry type
and are integral for integral~$k$, with $g_1=1$.
The parameter $a_k$ depends on the family in a natural way and is 
similar to the case of the $\zeta$--function.
The parameter $A$ depends both on the symmetry type
and the functional equation satisfied by the elements in the
family (specifically, it depends on the degree of the functional
equation in the relevant parameter).  Examples are given below.

One of the symmetry types described by Katz and Sarnak
is denoted O, for ``orthogonal.''  Examples of families
with this symmetry type conjecturally include

\itemitem{a)} the $L$--functions $L_f(s)$ associated 
with cusp forms $f\in S_m(\Gamma_0(1))$ of weight $m$ for the 
full modular group, 
\itemitem{b)} the $L$--functions $L_f(s)$ with $f\in S_2(\Gamma_0(N))$ 
of weight 2 
for the Hecke congruence group $\Gamma_0(N)$, 
\itemitem{c)} the twisted $L$--functions
$L(s,\hbox{sym}^{2\ell +1}(f)\otimes\chi_d)$ where 
$f$ is a self--dual cuspidal automorphic form on 
$GL_2$ and $\chi_d$ is a quadratic Dirichlet character mod~$|d|$, 
provided that the Fourier coefficients of $f$ have a Sato--Tate distribution, and  
\itemitem{d)} the twisted $L$--functions
$L_f(s,\chi_d)$ where $f$ is a self--dual cuspidal automorphic form on 
$GL_m$ for some $m$,
provided that the symmetric square $L$--function of $f$ 
does not have a pole at $s=1$.

To illustrate how the symmetry type is related to the mean values of
the $L$--functions in the family, we note the following conjectures.
$$\frac 1{m^*} \sum_{f\in S_m(\Gamma_0(1))}  L_f(
{\textstyle{\frac 12 }}
)^{k}\sim 
 g_k \frac{a_k}{\Gamma(1+ \frac 12 k(k-1))}
(\log m^\frac 12 )^{\frac 12 k(k-1)}$$
where $m^*$ is the cardinality of $S_m$, and
$$\frac 1{N^*} \sum_{f\in S_2(\Gamma_0(N))}  L_f(
{\textstyle{\frac 12 }}
)^{k}\sim 
 g_k \frac{a_k}{\Gamma(1+ \frac 12 k(k-1))}
(\log N^\frac 12 )^{\frac 12 k(k-1)}$$
where $N^*$ is the cardinality of $S_2(\Gamma_0(N))$, and also
$$\frac 1 {D^*}\sum_{|d|\le D} L(
{\textstyle{\frac 12 }}
,
\hbox{sym}^{2\ell+1}(f)\otimes\chi_d
)^{k}\sim 
 g_k \frac{a_k}{\Gamma(1+ \frac 12 k(k-1))}
(\log D^\frac 12 )^{\frac 12 k(k-1)},
$$
where $D^*$ is the number of quadratic characters
with conductor not exceeding $D$.  

Note:  these families can be further broken into the even forms and 
the odd forms, each of which is approximately half of the family. 
The average over the odd forms is identically zero
because the associated $L$--functions vanish at the center of
the critical strip.  

We were unable to locate a reference for these conjectures.
The above formulas were found, for example, using
the Petersson formula, and are based on results
found in \cite{D}\cite{DFI} and \cite{KMV}.  From those papers
we obtain $g_1=1$, $g_2=2$,
$g_3=2^3$ and $g_4=2^{7}$.

In the notation of (1.3), the families with symmetry type O
have $V(z)=z$ 
and $B(k)=\frac 12 k(k-1)$.  For the remaining parameter
we have $A=\Cal A$, where $\Cal A$ is the degree to which
the parameter $\Cal Q$ occurs in the functional equation.
For example, the first two families above satisfy a functional
equation of the form
$$
\Phi(s)=\left(\frac{N^{\frac 12}}{2\pi}\right)^{\!s}
 \Gamma\left(s+\frac{m-1}{2}\right) 
 L_f\left( s\right) 
=\varepsilon \bar{\Phi}(1-s) .
$$
That functional equation has degree $\frac 12$ in both
$N$ and $m$ aspect, so in the corresponding mean value we
have $A=\frac12$.

Another symmetry type considered by Katz and Sarnak is denoted Sp,
for ``Symplectic.''  Examples conjecturally are:
\itemitem{e)} Dirichlet $L$--functions $L(s,\chi_d)$, 
where $\chi_d$ is a quadratic Dirichlet character mod~$|d|$, 
\itemitem{f)} the symmetric square $L$--functions $L(s,
\hbox{sym}^2(f))$ associated 
with $f\in S_m(\Gamma_0(1))$,
\itemitem{g)} the twisted $L$--functions
$L(s,\hbox{sym}^{2\ell}(f)\otimes\chi_d)$ where 
$f$ is a self--dual cuspidal automorphic form on 
$GL_2$,
provided that the Fourier coefficients 
of $f$ have a Sato--Tate distribution, and  
\itemitem{h)} the twisted $L$--functions
$L_f(s,\chi_d)$ where $f$ is a self-dual cuspidal automorphic form on 
$GL_m$ for some $m$,
provided the symmetric square $L$--function of $f$ has a pole at $s=1$.

The conjectured mean value in case e) is:
$$\frac 1 {D^*}\sum_{|d|\le D} L(
{\textstyle{\frac 12 }}
,\chi_d)^{k}\sim 
 g_k \frac{a_k}{\Gamma(1+ \frac12 k (k+1))}
(\log D^\frac 12)^{\frac 12 k(k+1)},
$$
where $D^*$ is the number of quadratic characters
with conductor not exceeding $D$.  
This conjecture is based
on work of Jutila~\cite{J} and Soundararajan~\cite{S}.
In this case we have $g_1=1$, $g_2=2$,
$g_3=2^4$, and it is conjectured that
$g_4=3\cdot2^{8}$.  

In case f) we have the conjecture
$$\frac 1{m^*} 
\sum_{f\in S_m(\Gamma_0(1))}  L(
{\textstyle{\frac 12 }}
,\hbox{sym}^2(f))^{k}\sim 
 g_k \frac{a_k}{\Gamma(1+ \frac12 k (k+1))}
(\log m^{\frac 12})^{\frac 12 k(k+1)},
$$
where $m^*$ is the cardinality of $S_m$.
The cases g) and h) look just like e) above.

In the notation of (1.3), families with symmetry type Sp have $V(z)=z$ 
and $B(k)=\frac 12 k(k+1)$.  The parameter
$A$ for the Sp families is determined in exactly
the same way as for the O families.

A third symmetry type described by Katz and Sarnak is
denoted U, for ``Unitary.''
An example of a family with this symmetry type is
\itemitem{i)} $L(s,\chi)$ for $\chi$ a character mod~$q$

The conjectured mean value in this case is
$$\frac 1 {Q^*}\sum_{|q|\le Q} {\sum_{\chi\mod q}} |L(
{\textstyle{\frac 12 }}
,\chi)|^{2k}\sim  
g_k
\frac {a_k}{\Gamma(1+k^2)}
(\log Q)^{k^2},
$$
where the inner sum is over characters mod~$q$,
and $Q^*$ is the number of characters with 
conductor at most~$Q$.

For families with symmetry type U we have $V(z)=|z|^2$ 
and $B(k)=k^2$.  For the remaining parameter
we have $A=2\Cal A$, where $\Cal A$ is the degree to
which the parameter $\Cal Q$ appears in the functional equation. 
As the above example suggests, we can think of the 
Riemann $\zeta$--function as forming its own
Unitary family, where the mean values of the 
$\zeta$--function correspond to averages of special
values of the family $\{\zeta(\frac12 + it)\}_{t\in \R}$.

We end this discussion by giving an example of the
parameter $a_k$ in the above formulas.
The following is the $a_k$ associated to with
family of 
real Dirichlet $L$--functions $L(s,\chi_d)$,
which has symmetry type~Sp:
$$
a_k=\prod_p\frac{\left(1-\tfrac 1p\right)^{\tfrac{k(k+1)}{2}}}
{\left(1+\tfrac
1p\right)}
\Biggl(\frac12\({\left(1+\tfrac{1}{\sqrt{p}}\right)^{-k}+
\left(1-\tfrac{1}{\sqrt{p}}\right)^{-k}}\)+\frac 1p\Biggr) .
$$

While all of the above mean values are conjectural (except for certain
small values of~$k$), the various parameters in the formulas are
fairly well understood, except for the constant~$g_k$.  
Recently, Keating and Snaith \cite{KeSn} have used
techniques from random matrix theory to obtain conjectures
for $g_k$ for the symmetry types given above.  
This is described in the next section.

In this paper we study 
the functions $g_k$ 
in detail. 
We also report on calculations
of the mean square of the $L$--function times a
``mollifier'' for each of the symmetry types described above.
The results on $g_k$ are summarized in the next section,
and the results on mollifiers are presented in the following section.
The remainder of the paper is devoted to results related
to the various properties of the $g_k$ functions.

We thank J.~Keating, Z.~Rudnick, P.~Sarnak, A.~Selberg, and N.~Snaith for
helpful conversations.

\head
2. Statement of Results
\endhead

The constants~$g_k$
occurring in the conjectured mean values of these 
$L$--functions is the most mysterious aspect
of these mean values.
Conjectures for $g_k$ have recently been given, which we now describe.

The methods used to obtain 
conjectures for~$g_k$
are based on random matrix theory.
This approach to the study of $L$--functions
began with Montgomery's work \cite{M}
on the pair correlation of zeros of the Riemann $\zeta$--function,
and has been fruitful for establishing both rigorous and
conjectural results about $L$--functions.
The underlying idea is the 
conjecture that the zeros of $L$--functions 
are distributed on the critical line like the eigenvalues of 
matrices from the
Gaussian Unitary Ensemble of large random Hermitian matrices.
This is referred to as the ``GUE'' conjecture. 
This conjecture has been corroborated, to a large extent, 
by the extensive computations of Odlyzko \cite{O}. 
Heuristic explanations for the GUE conjecture
have been given by Bogolmony and Keating \cite{BK}. 
Now
Keating and Snaith \cite{KeSn} have announced conjectures for
$g_k$ that have been obtained by techniques from 
this theory.

Keating and Snaith's conjecture in the Unitary case is that 
$$
g_{\lambda,U} =\Gamma(1+B_U(\lambda))\lim_{N\to \infty} 
N^{-\lambda^2}
\prod_{j=1}^N 
\frac{\Gamma(j)\Gamma(j+2\lambda)}{\Gamma(j+\lambda)^2}.
\eqno(2.1)
$$
where $B_U(\lambda)=\lambda^2$.
The conjecture in the Orthogonal case is 
$$
g_{\lambda,O} =\Gamma(1+B_O(\lambda))\lim_{N\to \infty}
2^{2N\lambda} \prod_{j=1}^{N} \frac{\Gamma(N+j-1)
\Gamma(j-\frac12+\lambda)} 
{
\Gamma(N+j-1+\lambda)
\Gamma(j-\frac12) 
}
\eqno(2.2)
$$
where $B_O(\lambda)=\frac12 \lambda(\lambda-1)$.
And the conjecture in the Symplectic case is 
$$
g_{\lambda,Sp} =\Gamma(1+B_{Sp}(\lambda))
\lim_{N\to \infty}
2^{2N\lambda}
\prod_{j=1}^{N} \frac{\Gamma(N+j+1)\Gamma(j+\frac12+\lambda)
}  { 
\Gamma(N+j+1+\lambda)
\Gamma(j+\frac12) 
},
\eqno(2.3)
$$
with $B_{Sp}(\lambda)=\frac12 \lambda(\lambda+1)$.

We have recently learned that Br\'ezin and Hikami \cite{BH} have independently
obtained the above conjectures in the Sp and O cases.

It is of key importance that the above formulas agree with the known and 
conjectural values of $g_k$ given earlier.
These formulas for $g_\lambda$ suggest that it would be natural to
absorb the factor $\Gamma(1+B(\lambda))$ in to the definition
of~$g_\lambda$.  However, doing so would obscure the arithmetic 
origins of the parameters in the mean values, and we would lose
the significant fact that $g_k$ is integral for $k$ a
positive integer.

In this paper, we develop the properties of the 
Keating--Snaith constants $g_\lambda$.
We show that $g_\lambda/\Gamma(1+B(\lambda))$ 
is a nonvanishing  meromorphic function of order 2. 
We give a complete description of its pole locations,
and we express $g_\lambda/\Gamma(1+B(\lambda))$
in terms of the 
Barnes--Vign\'eras double $\Gamma$--function \cite{V}\cite{Sa}\cite{UN}.

Moreover, $g_k$ is an integer for positive integer $k$, 
and we establish asymptotic formulas for
$g_k$ as $k\to \infty$.
We remark that that there are interesting patterns  in the
prime factorization of~$g_{k}$.
The exponents of small primes dividing $g_{k}$ are quite irregular,
and in fact $p\nmid g_{k,U}$ for $k<p<k+\sqrt{p}$.
For example,  
$$\eqalign{ g_{100,U}=&\,2^{95}\cdot 3^{65}\cdot5^{24}\cdot7^{33}
\cdot 11^{10}\cdot13^{33}\cdot17^{36}\cdot19^{29}\cdot23^{20}
\cdot29^{16}\cdot31^{11}\cdot37^{10}\cdot41^{12}
\cr &\quad 
\cdot43^9\cdot47^4
\cdot53^3\cdot59^7\cdot61^9\cdot
67^{18}\cdot71^{12}\cdot 73^{10}\cdot 79^6\cdot83^4\cdot89^2
\cdot97\cdot113
\cr & \quad 
\cdot 127^5
\cdot131^7\cdot137^9
\cdot139^{10}\cdot149^{16}\cdot151^{17}\cdot
157^{20}\cdot 163^{24}\cdot167^{26}\cdot173^{30}
\cr &\quad 
\cdot179^{34}
\cdot 181^{36}\cdot191^{43}\cdot193^{44}
\cdot197^{47}
\cdot 199^{47}\cdot 211^{47}\cdot223^{44}\cdot \dots \cdot9973.
}$$

Even more interesting are the patterns in the exponent
of $p$ in the prime factorization of $g_k$ as $k\to\infty$.
We show that
there are continuous self--similar
functions $c_p(x)$ such that if
$k_j =[p^j x]$ then $v_p(g_{k_j,U})\sim k_j c_p(x)$, 
and $v_p(g_{k_j,O})\sim  v_p(g_{k_j,Sp})\sim \frac 12 k_j c_p(x)$,
where
$v_p(n)$ is the power of $p$ in the factorization of~$n$.
A graph of $c_3(x)$ is shown here:   
 
\vskip 0.2in

\hskip .4in \epsffile{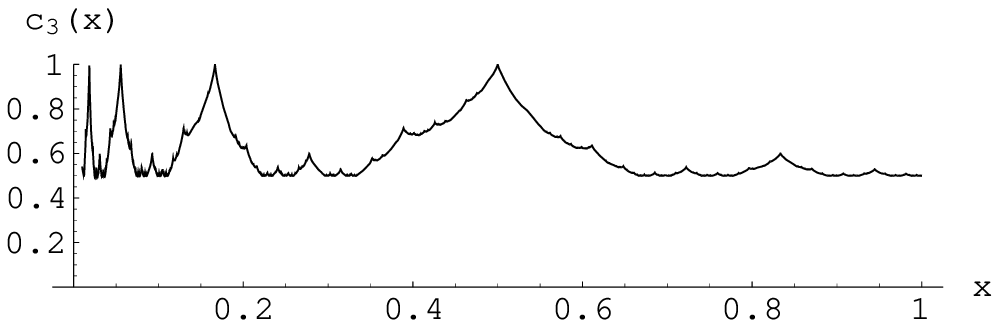}

\vskip 0.1in

Note that these functions satisfy $c_p(x)=c_p(px)$.
We will establish other self--similarity
properties of~$c_p(x)$, and also give an elegant
formula for~$c_p(x)$.

As a matter of interest, we also report that the Keating--Snaith
constant $g_{\frac 12,U }$ can be used to conjecture the mean 1st moment of
the $\zeta$--function, and a calculation finds:
$$\eqalignno{
g_{\frac 12,U }  & = \frac{\Gamma(\frac 54 )\pi^{\frac 14 }}
                        {2^{\frac 16 }}
\exp\left(
\frac 14\( \frac{\zeta'}{\zeta}(2)-\gamma +1\)
\right) 
&(2.4)\cr 
&\approx
1.0362329154....
}$$
It has been shown \cite{CG1}\cite{H--B} that 
$g_{\frac 12,U}$, if it exists,
satisfies
$1\le g_{\frac 12,U}\le 16/15$, assuming RH.

The paper is organized as follows.
In the next section we discuss the mean square of $L$--functions times
a mollifier.
The remainder of the paper is devoted to establishing properties
of $g_\lambda$ as defined in (2.1)--(2.3).  We will write
$g_\lambda$ and $B(\lambda)$ to mean $g_{\lambda,X}$ and $B_X(\lambda)$
for $X$ any one of U, O, or~Sp.  Results which are particular to 
one of the $g_{\lambda,X}$ will be specifically stated as such.
In section 4 we determine analytic properties 
of $g_\lambda$ as
a function of the complex number~$\lambda$.
This permits us to express $g_\lambda$ in terms of the
double $\Gamma$--function.
In section~5 we prove integrality properties of
$g_k$ for $k$ a positive integer.  
In section~6 we discuss the self--similar functions~$c_p(x)$.
In section~7 we 
give an asymptotic formula for $g_k$ as~$k\to\infty$.

\head
3. Mollified mean squares
\endhead

We model our discussion on  the mollified mean square 
required for Levinson's method \cite{Lev}\cite{C}.
If 
$$
L(s)^{-1}=\sum_{n=1}^\infty \frac{m(n)}{n^s},
$$
then our mollifier for $L(s)$ will be
$$
M(s,f)=\sum_{n\le y}\frac{m(n)}{n^s}
f\({\frac{\log y/n}{\log y}}\),
$$
where $f(x)$ is a real polynomial with $f(0)=0$,
and $y=Y^\theta$ for some $\theta>0$, where
$Y$ will be chosen appropriately for each family.
We are interested in an asymptotic formula for
the mean square of
$L(s)M(s,f)$ for $L$ in the families described in 
a previous section.

We first consider the Unitary family, for which our model is
the Riemann $\zeta$--function.
Let
$$
{\Cal M}_U(P,Q,\theta):=
\frac 1T \intl_1^T
\left|
Q\(\frac{-1}{\log T}\frac{d}{da}\)
\zeta(\thalf+a+it)
M(\thalf+it,P)
\right|^2_{a=0} dt .
$$
The length of the mollifier is~$T^\theta$.
The following formula holds \cite{C} for $\theta<\frac 47$:
$$\eqalign{
{\Cal M}_U(P,Q,\theta)
\sim &\,
P(1)^2
Q(0)^2 
+
\frac 1{\theta}
\intl_0^1\intl_0^1
\(
P'(x)
Q(y)
+\theta 
P(x)
Q'(y)
\)^2dx\,dy .
}$$
An interesting simple case is 
$\displaystyle{{\Cal M}_U(x,1,\theta)\sim 1+\theta^{-1}}$.

For the Orthogonal family, the mollifier has length 
$\sqrt{X}^{\,\theta}$.
The mollified mean value is:
$$
{\Cal M}_O(P,Q,\theta)=
\frac{1}{X^*} \sum_{c(f)\le X} \left.\left(Q\left(\frac 2{\log X}\frac
{d}{da}\right) \xi_f(1/2+a)M(1/2,P)\right)^2\right|_{a=0}.
$$
For example applications of related mollified mean values,
see \cite{IS}\cite{KM}\cite{KMV}. 
If $Q$ is even or odd, then for $\theta<1$ we have
$$
\eqalign{
{\Cal M}_O(P,Q,\theta)\sim &\,
\(P(1)Q'(1)+\frac1 {\theta}P'(1)Q(1)\)^2 \cr
&+
 \frac{1}{\theta}\intl_0^1\intl_0^1\(\frac1\theta P''(x) Q(y)-\theta
P(x)Q''(y)\)^2 dx\, dy .
}$$
In this case we have 
$\displaystyle{{\Cal M}_{O}(x,1,\theta)\sim \theta^{-2}}$.
It seems unexpected that $\displaystyle{{\Cal M}_{O}(x,1,\theta)\not\to1}$
as $\theta\to 0$.  This may be related to the distribution of
low--lying zeros of these $L$--functions~\cite{ILS}.

For the Symplectic family, the mollifier has length 
${\sqrt{X}}^{\,\theta}$.
The mollified mean value is:
$$
{\Cal M}_{Sp}(P,Q,\theta):=
\frac 1{X^*} \sum_{c(f)\le X} \left.\left(Q\left(\frac 2{\log X}\frac
{d}{da}\right) \xi_f(1/2+a)M(1/2,P)\right)^2\right|_{a=0} .
$$
Applications of related mean values can be found in~\cite{S}.
If $Q$ is odd then ${\Cal M}_{Sp}(P,Q,\theta)=0$.
If $Q$ is even, then for $\theta<1$ we have
$$\eqalign{
{\Cal M}_{Sp}(P,Q,\theta)
\sim&\,
\(P(1)Q(1) 
+\frac1 {\theta}P'(1) {\hat Q}(1)\)^2\cr
&
+
\frac{1}{\theta}\intl_0^1\intl_0^1\(\frac1\theta P''(x){\hat Q}(y)-\theta
P(x)Q'(y)\)^2 dx \, dy ,
}$$
where ${\hat Q}(y)=\int_0^yQ(u)~du$.
We have
$\displaystyle{{\Cal M}_{Sp}(x,1,\theta)\sim
(1+\theta^{-1})^2 }$.

We end this discussion by pointing out the beautiful relationship
$$
{\Cal M}_{Sp}(P,Q',\theta)
\sim
{\Cal M}_O(P,Q,\theta) ,
$$
which is transparent in the above formulas.
It was noted in \cite{S} that this relationship is the
source of the amazing ``coincidence'' 
of main results in the papers
\cite{S} and~\cite{KM}.
The same ``coincidence'' appears in the main result
of~\cite{CGG}, and it is plausible that there may be
a similar explanation for
that connection also.

\head
4. Meromorphicity of $g_\lambda$
\endhead
We prove 

\nobreak

\proclaim{Theorem 4.1} The function $g_\lambda/\Gamma(1+B(\lambda))$ is 
meromorphic of order $2$ in the whole complex plane. It never vanishes.
Furthermore, for $k=1,2,\dots,$ 
\item{}{U)}
$\ \ \displaystyle{\mathstrut g_{\lambda,U}\over \Gamma(1+B_U(\lambda))}$ 
has a pole of order $2k-1$ at
$\lambda=\frac12 -k$, and no other poles,
\item{}{O)} 
$\ \ \displaystyle{\mathstrut g_{\lambda,O}\over \Gamma(1+B_O(\lambda))}$ 
has a pole of order $k$ at
$\lambda=\frac12 -k$, and no other poles,
\item{}{Sp)} 
$\ \ \displaystyle{\mathstrut g_{\lambda,Sp}\over \Gamma(1+B_{Sp}(\lambda))}$ 
has a pole of order $k-1$ at
$\lambda=\frac12 -k$, and no other poles.

\endproclaim

The zero and pole locations of $g_\lambda/\Gamma(1+B(\lambda))$
indicate a connection with the $\Gamma$ and $\Gamma_2$--function.  
By combining Theorem~4.1 with the asymptotic 
formulas for $g_k$ given in Theorem~7.1, we obtain the following

\proclaim{Corollary 4.2} We have the following representations of 
$g_\lambda$
in terms of the $\Gamma$--function and the
Barnes--Vign\'eras double $\Gamma$--function:
\item{}
\item{}{U)} 
$\ \ \displaystyle{\mathstrut 
{g_{\lambda,U}\over \Gamma(1+B_U(\lambda))}
={2^{\frac1{12}}\,e^{3\zeta'(-1)}\,
e^{-2\lambda\zeta'(0)}\,
2^{-2\lambda^2}}
\,\frac
{\Gamma_2(\lambda+\frac12)^2}
{\Gamma(\lambda+\frac12)}
}, $
\item{}
\item{}
\item{}{O)} 
$\ \ \displaystyle{\mathstrut 
{g_{\lambda,O}\over \Gamma(1+B_{O}(\lambda))}
={2^{-\frac{17}{24}}\,e^{\frac32 \zeta'(-1) +\frac 12 \zeta'(0)}\,
e^{-\lambda\zeta'(0)}\, 2^{\lambda} \,
2^{-\frac 12 \lambda^2}}
\,{\Gamma_2(\lambda+\frac12)}},$
\item{}
\item{}
\item{}{Sp)} 
$\ \ \displaystyle{\mathstrut 
{g_{\lambda,Sp}\over \Gamma(1+B_{Sp}(\lambda))}
={2^{-\frac5{24}}\,e^{\frac32 \zeta'(-1) -\frac 12 \zeta'(0)}\,
e^{-\lambda\zeta'(0)}\, 2^{-\lambda}\,
2^{-\frac 12 \lambda^2}}
\,\frac
{\Gamma_2(\lambda+\frac12)}
{\Gamma(\lambda+\frac12)}
}.$
\item{}\mathstrut

In particular,
$g_{\lambda+1,O}=2^\lambda g_{\lambda,Sp}$, and
$$
{g_{\lambda,O}\over \Gamma(1+B_{O}(\lambda))}
{g_{\lambda,Sp}\over \Gamma(1+B_{Sp}(\lambda))}
=2^{\lambda^2-1}{g_{\lambda,U}\over \Gamma(1+B_U(\lambda))} .
$$

\endproclaim

See  \cite{V}\cite{UN} for details about the double $\Gamma$--function.

\demo{Proof of Theorem 4.1} We illustrate with the case of
$g_{\lambda,U}$, the other cases being similar.
Choose $J$ such that 
$1+|2\lambda|<J< |4\lambda|$. 
Then, for $j\ge J$ the real part of $j+2\lambda$ is positive. We will show that
$$
\lim_{N\to \infty} 
\biggl(
-\lambda^2\log N+\sum_{j=J}^N \log \Gamma(j)-2\log \Gamma(j+\lambda)+
\log \Gamma(j+2\lambda) 
\biggr)
\eqno (4.1)
$$
exists and is bounded by 
$|\lambda|^{2+\epsilon}$ for any $\epsilon >0$. 
This assertion, together with well--known properties of the 
$\Gamma$--function, imply that $g_{\lambda,U}/\Gamma(1+B_U(\lambda))$ 
is meromorphic of 
order at most~2,
with zeros and poles as described. 
It follows that the order of $g_{\lambda,U}/\Gamma(1+B_U(\lambda))$ 
is exactly~2 because the series
$
\sum{\rho^{-r}},
$
where $\rho$ runs through the set of poles, 
with multiplicity, has exponent of convergence~$r=2$.

Now we prove the above assertion.
We need to estimate
$$
\sum_{j=J}^N f_j(0)-2f_j(\lambda)+f_j(2\lambda) ,
$$
where 
$
f_j(x)=\log \Gamma(x+j).
$
By a mean value theorem,
$$
f_j(0)-2f_j(\lambda)+f_j(2\lambda)=\lambda^2 f_j''(\xi_j)
$$
for some $\xi_j$ between 0 and  $2\lambda$.

To evaluate 
$f_j''(\xi_j)$,
use the following well--known formula (see \cite{R}), 
valid for $s \ne 0$ and not on the negative real axis,
$$
\log\Gamma(s)=(s-1/2)\log s-s+\frac12\log 2\pi
+\frac 12 \intl_0^\infty\frac{\psi_2(u)+\frac16}{(u+s)^2}~du 
$$
with 
$$
\psi_2(u)=\{u\}^2-\{u\}
$$
and $\{u\}=u-[u]$.
We have
$$\eqalign{
f_j''(\xi)&=\frac{1}{\xi+j}-
\frac{1}{2(\xi+j)^2}-
6\intl_0^\infty\frac{\psi_2(u)}{(u+j+\xi)^4}~du\cr
&={1\over j}+O\({\xi/ j^2}\) .
}$$

It follows that the limit in (4.1) exists 
and is bounded by 
$\lambda^2\log J+O(\lambda/J)\ll 1+\lambda^{2+\varepsilon}$.
All the assertions are proven.
\enddemo

\head 
5. Integrality results
\endhead

In this section $k$
is a positive integer and we will be interested 
in integrality properties of~$g_k$.
By virtue of Corollary~4.2, these can be used to deduce 
rationality results for the double $\Gamma$--function.
Denote by $v_p(n)$ the power of $p$ dividing~$n$.

\proclaim{Theorem 5.1} For all prime~$p$ we have
$v_p(g_k)\ge 0$, so $g_k$ is an integer.  If
$p>B(k)$ then $v_p(g_k)=0$.  If $p<B(k)$ then
$$v_p(g_{k})=0
\hbox{\  if  and  only  if $p^2>B(k)$ and \ } 
\cases
k<p<k+\sqrt{p} &\ \ \ \ (U)\cr
k-\sqrt{\mathstrut k+p}<p<k+\sqrt{\mathstrut k+p} &\ \ \ \ (O),\cr
\endcases
$$
with the Sp case following from the O case and the 
relationship $g_{k+1,O}=2^k g_{k,Sp}$.
\endproclaim

We have the following useful formulas for $g_k$.

\proclaim{Lemma 5.2} For $k$ a positive integer we have
$$\eqalign{
g_{k,U}=&\, B_U(k)!\ 
\frac{\(\prod_{j=1}^{k-1}j!\)^2}{\prod_{j=1}^{2k-1}j!} \cr 
&=\,B_U(k)!\ 2^{k-k^2} \prod_{j=1}^{k-1}
\frac{1}{(2j-1)!!\,(2j+1)!!}\cr
g_{k,O}=&\,
B_O(k)!\ 
2^{B_O(k)+k-1}
\prod_{j=1}^{k-1}\frac{j!}{2j!} \cr
=&\,B_O(k)!\
2^{k-1}
\prod_{j=1}^{k-1}\frac1{(2j-1)!!} 
\cr
g_{k,Sp}=&\,
B_{Sp}(k)!\ 
2^{B_{Sp}(k)}
\prod_{j=1}^{k}\frac{j!}{2j!}  
\cr
=&\,
B_{Sp}(k)!\ 
\prod_{j=1}^{k}\frac{1}{(2j-1)!!}  
.
}$$ 
\endproclaim
Note that relationships between $g_{k,U}$, $g_{k,O}$, and $g_{k,Sp}$, 
as given at the end of Corollary~4.2,
can also be obtained from the above formulas.

\demo{Proof} 
The proofs are similar, so we illustrate with the
case of $g_{k,U}$.
We begin with 
formula~(2.1),  first using the fact that
$k$ is an integer. 
$$\leqalignno{
\frac{g_{k,U}}{B_U(k)!}&= N^{-k^2}\lim_{N\to \infty} 
  \prod_{j=1}^{N}\prod_{m=1}^k\frac{j+2k-m}{j+k-m}\cr
&=  \lim_{N\to \infty} \prod_{m=1}^k N^{-k}\,
  \frac{\prod_{j=N-k+1}^N(j+2k-m)}{\prod_{j=1}^k(j+k-m)}\cr
&=\lim_{N\to \infty} \prod_{m=1}^k\prod_{j=1}^k
N^{-1}\,
\frac{N-j+2k-m}{j+k-m}\cr
&=\prod_{m=1}^k\prod_{j=1}^k\frac{1}{j+k-m}\cr
&=\prod_{m=0}^{k-1}\prod_{j=1}^k\frac 1{j+m}\cr
&=
\prod_{j=0}^{k-1}\frac{j!}{(j+k)!} . 
}$$
The last expression leads easily to
the Lemma.
\enddemo

From Lemma~5.2 it follows that 
$p\nmid g_{k}$
for the cases listed in Theorem~5.1.
To prove the remaining statements
in Theorem~5.1 we will find a precise expression for $v_p(g_{k})$.
We only consider the cases of $g_{k,U}$ and $g_{k,O}$,
the remaining case following immediately from the relationship
$g_{k+1,O}=2^k g_{k,Sp}$.  Also, we only consider primes $p>2$
because it is easy to show that $g_k$ is divisible by a large power
of~2.
The main tool we need is the following Lemma.

\proclaim{Lemma 5.3} We have
$$\leqalignno{
v_p\(\prod_{j=1}^J j!\)=&\,\sum_{\ell=1}^\infty
(J+1) 
{\left[{J\over p^\ell}\right]}
-
{p^\ell\over 2}
\(
{\left[{J\over p^\ell}\right]}^2
+
{\left[{J\over p^\ell}\right]}
\)
\cr
&&\hbox{and}\cr
v_p\(\prod_{j=1}^J (2j-1)!!\)=&\,\sum_{\ell=1}^\infty
\(J+\frac12\)
{\left[{2J-1\over p^\ell}\right]_2}
-
{p^\ell\over 2}
\left[{2J-1\over p^\ell}\right]_2^2 
,
}$$
where $[x]_2=[\thalf ([x]+1)]$.
\endproclaim

\demo{Proof}  Regrouping the products gives
$$
\prod_{j=1}^J j! =
\prod_{j=1}^J j^{J-j+1} 
\ \ \ \ \ \ \ \ \ \ \ \ \ \ 
\hbox{and}
\ \ \ \ \ \ \ \ \ \ \ \ \ \
\prod_{j=1}^J (2j-1)!! =
\prod_{j=1}^J (2j-1)^{J-j+1}.
$$
Thus,
$$
v_p\Biggl(\prod_{j=1}^J j!\Biggr)=\sum_{\ell=1}^\infty
\sum_{n=1}^{\left[{J\over p^\ell}\right]}
(J-n p^\ell + 1),
$$
and
$$
v_p\Biggl(\prod_{j=1}^J (2j-1)!!\Biggr)=\sum_{\ell=1}^\infty
\sum_{n=1}^{\left[{2J-1\over p^\ell}\right]_2}
J+\frac{p^\ell+1}{2} -np^\ell .
$$
Evaluating the inner sums gives the Lemma.
\enddemo

Combining Lemmas 5.2 and 5.3 we obtain the following expressions for
$v_p(g_{k})$.

\proclaim{Proposition 5.4}  We have
$$\leqalignno{
v_p(g_{k,U})= &\, \sum_{\ell=1}^\infty
\Biggl(
\[{k^2\over p^\ell}\]+(2k-p^\ell)\[{k-1\over p^\ell}\]
+\({p^\ell\over 2}-2k\)\[{2k-1\over p^\ell}\] 
\cr
&\phantom{xxxxx}
-p^\ell\[{k-1\over p^\ell}\]^2
+{p^\ell\over 2}\[{2k-1\over p^\ell}\]^2 
\Biggr)\cr
&&\hbox{and}\cr
v_p(g_{k,O})= &\, \sum_{\ell=1}^\infty
\[{\frac12 k(k-1)\over p^\ell}\] -
\(k-\frac12\)\[{2k-3\over p^\ell}\]_2+
\frac{p^\ell}{2}\[{2k-3\over p^\ell}\]_2^2  .
}$$
\endproclaim

To complete the proof we must obtain lower bounds for the 
above expressions.  The $g_{k,O}$ case is slightly simpler, 
so we handle it first.

Denote the summand in the second statement of Proposition~5.4 
by~$v_{p,\ell}(g_{k,O})$.  
Use $\theta$ to denote a number $0\le \theta < 1$.
For example, $[x]=x-\theta$.
We have:
$$\eqalignno{
v_{p,\ell}(g_{k,O})=&\,
\[{\frac12 k(k-1)\over p^\ell}\] -
\(k-\frac12\)\[{2k-3\over p^\ell}\]_2+
\frac{p^\ell}{2}\[{2k-3\over p^\ell}\]_2^2\cr
=&\, \frac1{2p^{\ell}}\(k^2 -  k\) -\theta -
\(k-\frac12\)\[{2k-3\over p^\ell}\]_2+
\frac{p^\ell}{2}\[{2k-3\over p^\ell}\]_2^2\cr
=&\, \frac1{2p^\ell}\(
\(k-p^\ell \[{2k-3\over p^\ell}\]_2\)^2-
\(k-p^\ell \[{2k-3\over p^\ell}\]_2\)\) -\theta 
&(5.1) \cr
>&\,-1,
}$$
because $M^2-M\ge0$ for integral $M$.  Thus, $v_{p,\ell}(g_{k,O})\ge 0$
for all $p$ and $\ell$, so $g_{k,O}$ is an integer.
The remaining statements about $v_{p}(g_{k,O})$ can be obtained
by considering the possibilities for $v_{p,\ell}(g_{k,O})$.

For the case of $g_{k,U}$,
denote the summand in the first statement of Proposition~5.4 
by~$v_{p,\ell}(g_{k,U})$.  
There are two cases:

\demo{Case 1}  $\displaystyle{\[{2k-1\over p^\ell}\]=
2\[{k-1\over p^\ell}\]}$.  We have
$$\eqalignno{
v_{p,\ell}(g_{k,U})&= 
\[{k^2\over p^\ell}\]
-2k\[{k-1\over p^\ell}\]
+p^\ell\[{k-1\over p^\ell}\]^2 \cr
&= {k^2\over p^\ell}-\theta
-2k\[{k-1\over p^\ell}\]
+p^\ell\[{k-1\over p^\ell}\]^2 \cr
&=p^{-\ell}\( k-p^\ell \[{k-1\over p^\ell}\]\)^2 -\theta .
&(5.2)\cr
&> -1 .
}$$
Thus, $v_{p,\ell}(g_{k,U})\ge 0$.
\enddemo

\smallskip

\demo{Case 2}  $\displaystyle{\[{2k-1\over p^\ell}\]=
2\[{k-1\over p^\ell}\]+1}$.  We have
$$\eqalignno{
v_{p,\ell}(g_{k,U})&= 
\[{k^2\over p^\ell}\]
+(2p^\ell-2k)\[{k-1\over p^\ell}\]
+p^\ell\[{k-1\over p^\ell}\]^2 +p^\ell - 2k\cr
&= {k^2\over p^\ell}-\theta
+(2p^\ell-2k)\[{k-1\over p^\ell}\]
+p^\ell\[{k-1\over p^\ell}\]^2 +p^\ell - 2k\cr
&=p^{-\ell}\(k-p^\ell - p^\ell \[{k-1\over p^\ell}\]\)^2 -\theta .
&(5.3)\cr
&> -1 .
}$$
Again, $v_{p,\ell}(g_{k,U})\ge 0$. 
\enddemo

This covers all cases because for integral $m$, $n$,
we have
$$
{{n\over m}-1+{1\over m}
\le \[{\mathstrut n\over\mathstrut m}\] \le {n\over m}+1-{1\over m}}.
$$

This proves that $g_{k,U}$ is an integer.  The remaining statements
can be proven by letting $\ell=1$ and considering
the possible cases.  This completes the proof of Theorem~5.1.

\head 
6. The function $c_p$ 
\endhead

In this section we study the prime factorization of
$g_k$ as $k\to\infty$.
We find that $v_p(g_{k})$ is described by an interesting 
self--similar function.

\proclaim{Theorem 6.1} 
Define
$$
c_p(x)=x^{-1}\sum_{\ell=-\infty}^\infty 
p^{-\ell} ||p^\ell x||^2 ,
$$
where $||x||$
is the distance from $x$ to the nearest integer.
Let $x>0$ and put $k_j = [p^j x]$.  Then as $j\to\infty$,
$$
v_p(g_{k_j,U})\sim k_j c_p(x)
\ \ \ \ \ \ \ \ \ \ \ \ 
\hbox{and}
\ \ \ \ \ \ \ \ \ \ \ \ 
v_p(g_{k_j,O})\sim v_p(g_{k_j,Sp})\sim\frac 12  k_j c_p(x) .
$$
\endproclaim
We will leave to the reader the exercise of showing
that $c_2(x)=1$.  For the remainder of this section
we will assume $p\ge 3$.
Here are some graphs:

\vskip 1.6in

\hskip -0.1in \epsffile{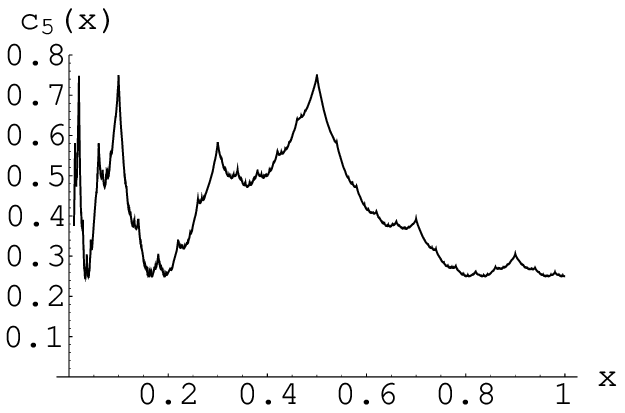}
\hskip .3in \epsffile{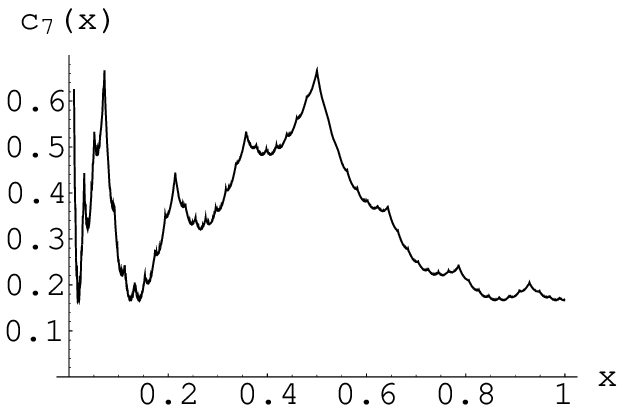}
\vskip -1.4in

\noindent Note that as $p\to\infty$, $c_p(x)$ 
approaches $x^{-1}||x||$, uniformly for
$\delta\le x\le \delta^{-1}$.

From the graphs it appears that each $c_p(x)$ is not differentiable.
Much more is actually true.  At each rational point $a/b$,
the function $c_p(x)$ is either self--similar, or
it has a cusp, or it has a vertical tangent.  This 
is made precise in Theorem~6.2.

Let $[[{n}]]$ denote the absolute least
residue of $n\bmod b$.  That is, $n\equiv[[{n}]]\bmod b$
and $b/2<[[{n}]]\le b/2$.

\proclaim{Theorem 6.2} Suppose $p$, $a$, and $b$ are pairwise
coprime, and suppose $p^r\equiv 1\bmod b$,
with $r>0$.  
Let $f(x)=||x||^2$.

If 
$\displaystyle{\ \sum_{j=0}^{r-1}[[{ap^j}]]=0}$, then
$c_p(x)$ has the following self--similarity property 
at $x=a/b:$
$$\eqalign{
\lim_{{m\to\infty}\atop {m\equiv m_0\bmod r}}
\frac{c_p\left( \frac ab +p^{-m}x \right) - c_p\left(\frac ab \right)}
{p^{-m}x}
= &\, 
\frac {x}{p-1} 
+ 
\sum_{\ell=1-m_0}^\infty f'\left(
p^{-\ell}\frac a{b}\right)
\cr
&+
x^{-1}
\sum_{\ell=0}^\infty                                                      
p^{-\ell}\left(                                                                 
f\left(p^{\ell+m_0}\frac  ab
+p^\ell x\right)-f\left(p^{\ell+m_0}\frac ab\right)
\right)   .
}$$
\nobreak
Note that this holds if\/ $b=1$, or more
generally if\/ $b\not=2$ and\/ $p^n\equiv -1\bmod b$ 
for some~$n$.

\vbox{
Suppose the above condition does not hold.
In the case $b=2$ we have
$$
c_p\left(\frac a2 + p^{-m}x \right) =
c_p\left( \frac a2 \right) - mp^{-m}|x| + O(p^{-m}x) ,
$$
so $c_p$ has a cusp at $a/2$.
For all other $b$ there is a nonzero constant $k$ so that
$$
c_p\left(\frac ab + p^{-m}x \right) =
c_p\left( \frac ab \right) + k m p^{-m}x + O(p^{-m}x) ,
$$
so $c_p$ has a vertical tangent at $a/b$.
}
\endproclaim
Note:  if $(p,b)\not=1$ then one can use the relation
$c_p(x)=c_p(p^j x)$ and then apply the Theorem at $p^j a/b$.

We present some graphs of $c_p(x)$ near rational points $a/b$ where they
are self--similar.  The simplest case is is when $b=1$,
for then $c_p(x)$ is self--similar under scaling by~$1/p$.
We give the examples of $c_3(x)$ near $x=8$ and 
$c_7(x)$ near $x=3$.
In each graph  the $x$--axis is interpreted as 
extending $1/p^m$ on either side of the central point,
for any large~$m$, and the vertical scale is also
on the order of~$1/p^m$.

\vbox{
\phantom{oo}
\vskip 2.5in
\hskip -0.1in \epsffile{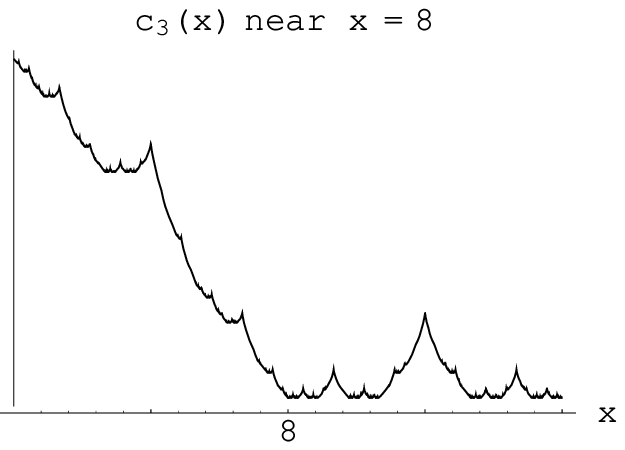}
\hskip .3in \epsffile{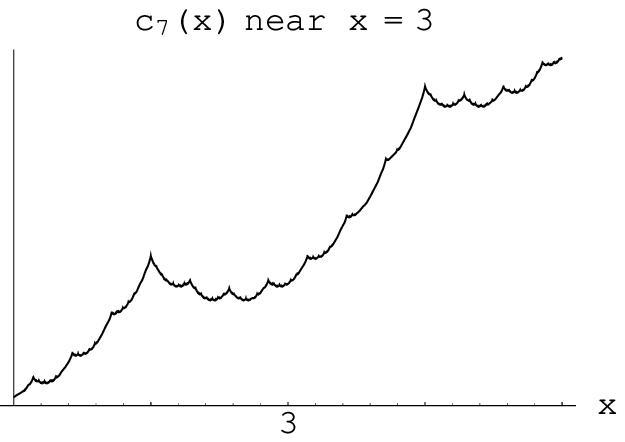}
\vskip -2.45in
}

Our last example is $c_5(x)$ near $x=3/13$.  Since $5$ has
order $4\bmod 13$, and $5^2\equiv -1\bmod 13$, 
the function is self--similar on
rescaling by~$1/p^4$.  Thus, we have four possible
pictures.  In the graphs below, the $x$--axis is 
interpreted as extending $1/p^m$ on either side of the
central point, with $m\equiv i\bmod 4$ for the 
graph in the $i$th quadrant.  The graphs should be 
read counterclockwise, with each successive graph
being the middle $1/5$th of the previous graph.

\nobreak

\vbox{
\phantom{oo}
\vskip -0.5 in
\hskip -0.2in \epsffile{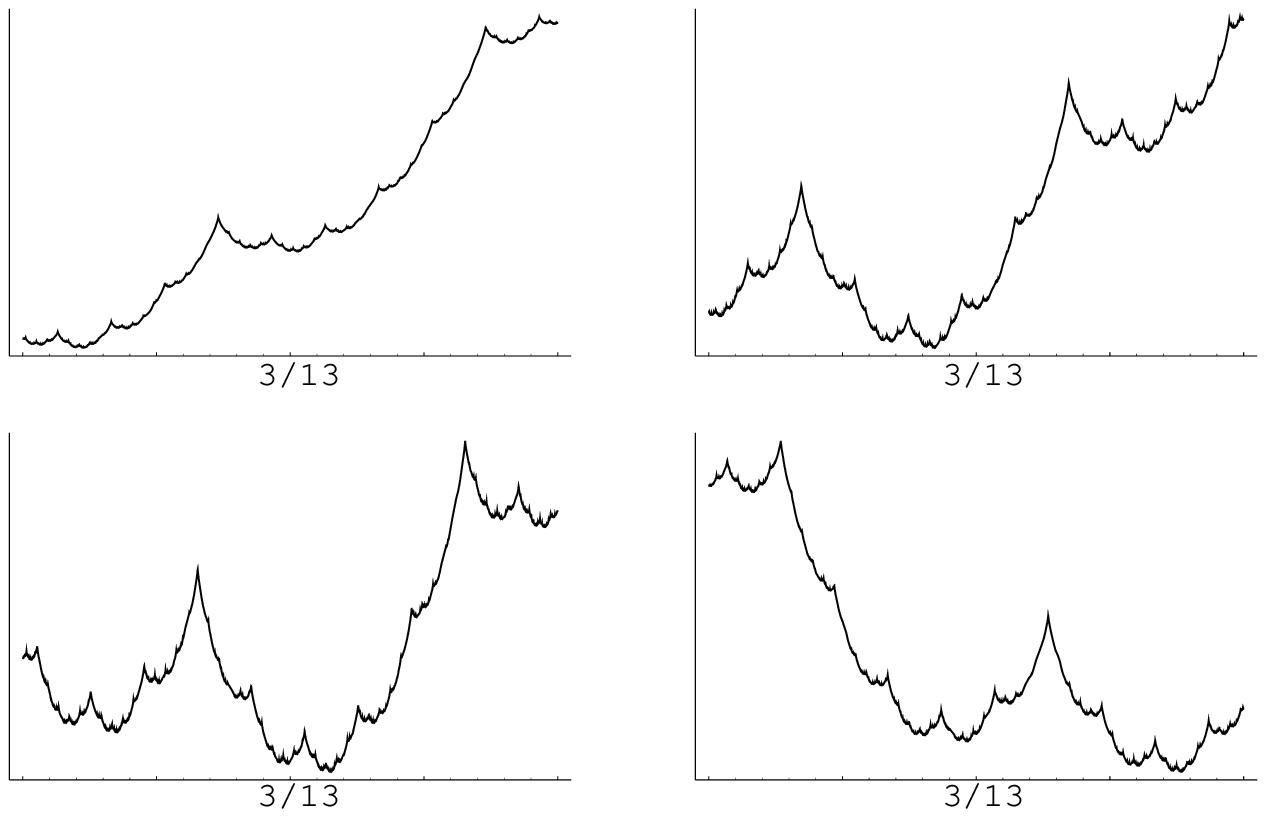}
\vskip 0.5in
}
Note that $c_5(3/13)=23/72$.  It is not difficult to evaluate 
$c_p(a/b)$ for any particular $a/b$, and in addition one
sees that $c_p(a/b)$ is rational.

The self--similarity properties of $c_p(x)$ 
follow from the formula for $c_p(x)$ in
Theorem~6.1.  We have

\proclaim{Proposition 6.3} Let $f$ be continuous on $\R$, periodic with period $1$,
and twice continuously differentiable except at finitely
many points $(\bmod\ 1)$ at which it is twice 
differentiable
from both the left and the right.
Also suppose $f(x)\ll x^2$ and
$f'(x)\ll |x|$
as $x\to 0$. Define
$$
d(x)=\sum_{\ell=-\infty}^\infty 
p^{-\ell}{f\left(p^\ell x\right)}.
$$
Then $d(x)$ is continuous on $\R$ and satisfies
$d(px)=pd(x)$.  
Also suppose
$p^r\equiv 1 \mod b$ and $m\equiv m_0\bmod r$,
with $m_0\ge 0$. If $x$ is sufficiently small then
$$\eqalign{
{d\left(\frac ab+ {p^{-m}}x\right)-d\left(\frac
ab\right)}
= &\, p^{-m}x\sum_{\ell=1-m_0}^\infty f_\pm'\left(
p^{-\ell}\frac a{b}\right)
\cr
& +
p^{-m}\sum_{\ell=0}^\infty
p^{-\ell}\left(
{f\left(p^{\ell+m_0}\frac
ab+p^\ell x\right)-f\left(p^{\ell+m_0}\frac ab\right)}
\right)
\cr
& +
(m-m_0)  p^{-m} x{1\over r}\sum_{j=0}^{r-1} 
f_\pm'\left(
p^{j}\frac ab \right)
\cr
& +
{{\frac 12}}  p^{-2m}x^2\sum_{\ell=0}^{m-1} p^{\ell} 
f_\pm''\left(
p^{\ell}\frac ab +p^{\ell-m}\xi\right) +O(p^{-2m}),
}$$
where $\xi=\xi_\ell$ is between $0$ and~$x$.
Here the $\pm$ refers to a left-- or right--derivative
and has the same sign as~$x$. 
\endproclaim

Note:  the condition that $x$ be sufficiently small can be
made explicit if one specifies the points at which $f$
is not smooth.

\demo{Proof of Theorem 6.2}
Apply Proposition~6.3 with $f(x)=||x||^2$.
We have
$$
f_\pm'\left(\frac ab \right)
=
\cases 
f'(\frac ab) &\hbox{ if }\   b\not=2 \cr
\mp 1 &\hbox{ if }\  b=2,
\endcases
$$
and $f_\pm''(\frac ab )=2$ for all $\frac ab $.
This makes the sum over $f_\pm''$ explicit.
For the second sum over $f_\pm'$ note that
$$
\sum_{j=0}^{r-1} [[{ap^j}]]=
\frac b2
\sum_{j=0}^{r-1} f'\left( p^j\frac ab \right) .
$$
This proves Theorem~6.2.
\enddemo

It remains to prove Theorem~6.1 and Proposition~6.3.

\demo{Proof of Proposition 6.3} 
The continuity of $d(x)$ follows from the Weierstrass $M$--test,
and the functional equation $d(p x)=pd(x)$ follows
by changing the summation index.

For the remaining properties,
write
$$
d\left(\frac ab \right) - 
d\left(\frac ab +p^{-m}x \right)   
=  A+B+C,
$$
where we used the definition of $d$ and
split the resulting sum into three pieces: $0\le \ell < m$, 
$m\le \ell < \infty$, and $\ell <0$. 

We have
$$\eqalign{
B = &\, 
\sum_{\ell=m}^{\infty} p^{-\ell} \left( f\left(
p^{\ell}\frac ab +p^{\ell-m}x \right)
-
f\left(
p^{\ell}\frac ab \right)\right) \cr
= &\, 
p^{-m}
\sum_{\ell=0}^{\infty} p^{-\ell} \left( f\left(
p^{\ell+m}\frac ab +p^{\ell}x \right)
-
f\left(
p^{\ell+m}\frac ab \right)\right) \cr
= &\, 
p^{-m}
\sum_{\ell=0}^{\infty} p^{-\ell} \left( f\left(
p^{\ell+m_0}\frac ab +p^{\ell}x \right)
-
f\left(
p^{\ell+m_0}\frac ab \right)\right)  .
\cr
}$$
We did a change of variable, then used the assumption
$p^r\equiv 1\bmod b$, the periodicity of $f$,
and the fact that $m_0, \,\ell\ge 0$.

Next,
$$\eqalign{
C = &\, 
\sum_{\ell=1}^{\infty} p^{\ell}\left( f\left(
p^{-\ell}\frac ab +p^{-\ell-m}x \right)
-
f\left(
p^{-\ell}\frac ab \right)\right) \cr
= &\, 
\sum_{\ell=1}^{\infty} p^{\ell}\left( 
f_\pm'\left(
p^{-\ell}\frac ab \right)
p^{-\ell-m}x
+O\left(
p^{-2\ell -2 m}
f_\pm''\left(
p^{-\ell}\frac ab \right)
\right)
\right)
\cr
= &\, 
p^{-m}
x
\sum_{\ell=1}^{\infty} 
f_\pm'\left(
p^{-\ell}\frac ab \right)
+O\left(
p^{-2 m}
\right)  .
}$$
We used Taylor's theorem and the fact that $f''$ is bounded.

Next,
$$\eqalign{
A = &\, 
\sum_{\ell=0}^{m-1} p^{-\ell} \left( 
f\left(
p^{\ell}\frac ab +p^{\ell-m}x \right)
-
f\left(
p^{\ell}\frac ab \right)
\right)
\cr 
= &\,
\sum_{\ell=0}^{m-1} p^{-\ell} \left( 
f_\pm'\left(
p^{\ell}\frac ab \right)
p^{\ell-m}x +
\frac 12 f_\pm''\left(
p^{\ell}\frac ab +p^{\ell-m}\xi\right) p^{2\ell-2m}x^2
\right)
\cr 
= &\,
 p^{-m}x\sum_{\ell=0}^{m-1} 
f_\pm'\left(
p^{\ell}\frac ab \right)
+
\frac 12 x^2 p^{-2m}\sum_{\ell=0}^{m-1} p^{\ell} 
f_\pm''\left(
p^{\ell}\frac ab +p^{\ell-m}\xi\right) ,
}$$
where $\xi=\xi_\ell$ is between $0$ and~$x$.  
This is valid if $x$ is sufficiently small, 
in terms of $a/b$, $p$, and~$f$.
Now,
$$\eqalign{
\sum_{\ell=0}^{m-1}
f_\pm'\left(
p^{\ell}\frac ab \right)
= &\,
\sum_{j=0}^{r-1}
\sum_{\ell=0}^{{{m-m_0}\over r}-1}
f_\pm'\left(
p^{r\ell+j}\frac ab \right)
+
\sum_{j=0}^{m_0-1}
f_\pm'\left(
p^{m-m_0+j}\frac ab \right)
\cr
= &\,
\sum_{j=0}^{r-1}
\sum_{\ell=0}^{{{m-m_0}\over r}-1}
f_\pm'\left(
p^{j}\frac ab \right)
+
\sum_{j=0}^{m_0-1}
f_\pm'\left(
p^{j}\frac ab \right)
\cr
= &\,
\frac {m-m_0}r 
\sum_{j=0}^{r-1}
f_\pm'\left(
p^{j}\frac ab \right) 
+
\sum_{j=0}^{m_0-1}
f_\pm'\left(
p^{j}\frac ab \right).
}$$
This completes the proof of Proposition~6.3.

\enddemo

\demo{Proof of Theorem 6.1}  
First we consider the Unitary case.
We may unify formulas (5.2) and (5.3) to have
$$
v_p(g_{k,U})=\sum_{\ell=1}^\infty
p^{-\ell}\(k-p^\ell\(
\[{2k-1\over p^\ell}\]-\[{k-1\over p^\ell}\]\)\)^2
+O\(\log k \).
$$
We have
$$\eqalignno{
\lim_{j\to\infty} {v_p(g_{[p^j x],U})\over p^j }  
& =\lim_{j\to\infty} 
{1\over p^j }\sum_{\ell=1}^\infty
p^{-\ell}\([p^j x]-p^\ell\(
\[{2[p^j x]-1\over p^\ell}\]-\[{[p^j x]-1\over p^\ell}\]\)\)^2
\cr
& =
\lim_{j\to\infty} 
{1\over p^j }\sum_{\ell=1}^\infty
p^{-\ell}\(p^j x -\theta_j -p^\ell\(
\[{2p^j x-2\theta_j-1\over p^\ell}\]-\[{p^j x-1\over p^\ell}\]\)\)^2
\cr
& =
\lim_{j\to\infty} 
{1\over p^j }\sum_{\ell=1}^\infty
p^{\ell}\Bigl(p^{j-\ell} x -
\[2p^{j-\ell} x-(2\theta_j+1)p^{-\ell}\]+\[p^{j-\ell} x-p^{-\ell}\]\Bigr)^2
\cr
& =
\lim_{j\to\infty} 
\sum_{\ell=-\infty}^{j-1}
p^{-\ell}\Bigl(p^{\ell} x -
\[2p^{\ell} x-(2\theta_j+1)p^{\ell-j}\]+\[p^{\ell} x-p^{\ell-j}\]\Bigr)^2 .
}$$
We first used $[p^jx]=p^jx -\theta_j$ with $0\le \theta_j<1$,
 and $\[[x]/n\]=[x/n]$.
We eliminated the first $\theta_j$ by multiplying out and verifying
that the other terms make no contribution.
Then we changed variables $\ell \mapsto j-\ell$

Let 
$$a_{p,j}(x)= 
\sum_{\ell=-\infty}^{j-1}
p^{-\ell}\Bigl(p^{\ell} x -
\[2p^{\ell} x-(2\theta_j+1)p^{\ell-j}\]+\[p^{\ell} x-p^{\ell-j}\]\Bigr)^2 
$$
denote the expression inside the above limit.

Let $\lfloor x\rfloor$ denote the largest integer strictly
smaller than~$x$.  Note that 
$$
[t-\delta]=\lfloor t \rfloor
\ \ \ \ \ \ \ \ \
\hbox{if}
\ \ \ \ \ \ 0<\delta<t-\lfloor t \rfloor .
$$
Also, there is always such a $\delta$ because
$t\not=\lfloor t\rfloor$ for all~$t$.
We will use this to simplify $a_{p,j}(x)$.

Let $\varepsilon>0$ be given and
choose $L>0$ so that $p^{-L}<\varepsilon$.
Then choose $\delta>0$ such that
$$
[2 p^\ell x -\delta]=\lfloor 2p^\ell x\rfloor
\ \ \ \ \ \ \ \ \ \
\hbox{and}
\ \ \ \ \ \ \ \ \ \ 
[p^\ell x -\delta]=\lfloor p^\ell x\rfloor
$$
for $1\le \ell \le L$.  Now choose $J$ so that
$3 p^{L-J}<\delta$.  If $j\ge J$ we have
$$\eqalign{
a_{p,j}(x) & = 
\sum_{\ell=-\infty}^{L}
p^{-\ell}\Bigl(p^{\ell} x -
\lfloor 2p^{\ell} x\rfloor +\lfloor p^{\ell} x\rfloor \Bigr)^2 
\cr
&\phantom{x}+
\sum_{\ell=L+1}^{j-1}
p^{-\ell}\Bigl(p^{\ell} x -
\[2p^{\ell} x-(2\theta_j+1)p^{\ell-j}\]+\[p^{\ell} x-p^{\ell-J}\]\Bigr)^2 
\cr & =
\sum_{\ell=-\infty}^{\infty}
p^{-\ell}\Bigl(p^{\ell} x -
\lfloor 2p^{\ell} x\rfloor +\lfloor p^{\ell} x\rfloor \Bigr)^2 
+O\(\sum_{\ell=L+1}^\infty p^{-\ell}
+\sum_{\ell=L+1}^{J-1}  p^{-\ell}\)\cr
 & =
\sum_{\ell=-\infty}^{\infty}
p^{-\ell}\Bigl(p^{\ell} x -
\lfloor 2p^{\ell} x\rfloor +\lfloor p^{\ell} x\rfloor \Bigr)^2 
+O(p^{-L}) ,
}$$
and we were given that $p^{-L}<\varepsilon$.
Thus,
$$\eqalign{
\lim_{j\to\infty} a_{p,j}(x)=&\,
\sum_{\ell=-\infty}^{\infty}
p^{-\ell}\Bigl(p^{\ell} x -
\lfloor 2p^{\ell} x\rfloor +\lfloor p^{\ell} x\rfloor \Bigr)^2 
.
}$$
To finish the proof, note that
$(x-\lfloor 2x\rfloor +\lfloor  x\rfloor )^2 
=(x-[2x]+[x])^2=||x||^2$.

For the Orthogonal case, a similar calculation shows
$$\eqalign{
\lim_{j\to\infty} {v_p(g_{[p^j x],O})\over p^j } =&\,
 \frac12 \sum_{\ell=-\infty}^{\infty}
p^{-\ell}\(p^{\ell} x -
[ 2p^{\ell} x]_2\)^2 \cr
=&\,\frac12\, x\, c_p(x),
}$$
since $(x-[2x]_2)^2=||x||^2 $.
The same holds in the Symplectic case because 
$g_{k+1,O}=2^k g_{k,Sp}$.  This completes the proof of
Theorem~6.1

\enddemo

\head
7. Asymptotics
\endhead 

In this section we determine the asymptotic behavior of 
$g_k$ for large (integer) $k$. 
We have

\proclaim{Theorem 7.1} As $k\to \infty$, 
$$\eqalign{ \log g_{k,U}=&\,\log B_U(k)!
-k^2\log k+\(\frac32-2\log2\)k^2
-\frac1{12}\log k \cr
&\,\qquad +\frac1{12}\log  2 +\zeta'(-1)+O\left({k^{-1}}\right)
\cr
=&\,k^2\log k+\(\frac12-2\log2\)k^2+\frac {11}{12} \log k
\cr
&\,\qquad 
+\frac1{12}\log 2 -\zeta'(0)+\zeta'(-1)+
O\left({k^{-1}}\right),\cr
\log g_{k,O}=&\,
\log B_O(k)!
-\frac 12 k^2\log k + \(\frac 34 - \frac 12 \log 2\)k^2
      + \frac 12 k\log k + \(\log 2 - \frac12\)k \cr
    &\, \qquad - \frac 1{24}\log k - \frac{17}{24}\log 2 + \frac 12 \zeta'(-1)
+O\left({k^{-1}}\right)
\cr
=&\,\frac12 k^2\log k + \(\frac14 - \log 2\)k^2 -
      \frac12 k\log k + \(\frac32 \log 2 - \frac12\)k +
      \frac{23}{24}\log k 
           \cr
    &\, \qquad  - \frac{29}{24}\log 2 + \frac 1{4} 
  - \zeta'(0) + \frac 12 \zeta'(-1) +O\left({k^{-1}}\right),\cr
\log g_{k,Sp}=&\,
\log B_{Sp}(k)!
-\frac12 k^2\log k + \(\frac34 - \frac12\log 2\)k^2
      - \frac12  k\log k + \(\frac12 - \log 2\) k \cr
&\, \qquad  
- \frac{1}{24}\log k - \frac5{24}\log 2 
  + \frac12 \zeta'(-1)+O\left({k^{-1}}\right)
\cr
=&\,
\frac12 k^2\log k + \(\frac14 - \log 2\)k^2 +
      \frac12 k\log k + \(\frac12 - \frac32 \log 2\)k 
+ \frac{23}{24}\log k \cr
    &\, \qquad 
- \frac{17}{24}\log 2 + \frac 14 
  - \zeta'(0) + \frac 12 \zeta'(-1)
+O\left({k^{-1}}\right).
}$$
\endproclaim

Combining these with the asymptotics of the $\Gamma$ and
$\Gamma_2$--function
\cite{V}\cite{UN}
we obtain the expressions for $g_{\lambda}$ in terms of
$\Gamma_2$ given in Corollary~4.2.

To prove Theorem~7.1, combine the expressions
for $g_\lambda$ given in Lemma~5.1 with the elementary
formulas
$$
\sum_{j=1}^n\log j!=
(n+1)\sum_{j=1}^n\log j
-\sum_{j=1}^n j \log j,
$$
and
$$\eqalign{
\sum_{j=1}^n\log 2j!=&\,
\frac12 n(n+1)\log 2+
(n+1)\sum_{j=1}^n\log j 
+(n+1)\sum_{j=1}^n \log(2j-1)\cr
&\,
-\sum_{j=1}^n j\log j 
-\sum_{j=1}^n j\log(2j-1) .
}$$
The formulas in the following Lemma are sufficient to complete
the calculations.

\proclaim{Lemma 7.2} As $n\to \infty$ we have
$$\eqalign{
\sum_{j=1}^n \log j =&\, n\log n-n +\frac12 \log n  -\zeta'(0)
+\frac1{12n}+O\left({n^{-2}}\right)\cr
\sum_{j=1}^n \log (2j-1)=&\,
n\log 2n-n +\frac12 \log 2 -\frac1{24n}+
O\left( {n^{-2}}\right)\cr
\sum_{j=1}^n j\log j
=&\,\frac12 n^2\log n - \frac 14 n^2
+\frac12 n\log n
+
\frac{1}{12}\log n
+ \frac 1{12}-\zeta'(-1) +O(n^{-1})\cr
\sum_{j=1}^n j\log (2j-1)
=&\,
\frac12 n^2\log 2n  - \frac 14 n^2
+\frac12 n\log 2n
 -\frac12 n
-
\frac{1}{24}\log n \cr
&+\frac{7}{24}\log 2
- \frac 1{24}+\frac12 \zeta'(-1) +O(n^{-1})
}$$
\endproclaim

The first formula in Lemma 7.2 is Stirling's formula,
and each expression can be derived from the Euler--Maclaurin 
summation formula
(see Rademacher \cite{R}).
This completes the proof of Theorem~7.1.

\enddocument